\documentclass[12pt]{article}
\usepackage{amssymb}
\newtheorem{lemma}{Lemma}
\newtheorem{theorem}{Theorem}

\title{Prime sum graphs and the induced trees they contain}
\author{Ernie Croot and Patrick Jin}

\begin{document}

\maketitle

\begin{abstract} In this paper we show that prime sum graphs on $n$ vertices -- which are graphs on vertex set $\{1,2,...,n\}$ where $ij$ is an edge when $i+j$ is prime -- contain all trees with at most $\exp( c \log n / \log\log n)$ vertices as induced subgraphs.  We also prove some results for related graphs, and end with some unsolved problems.
\end{abstract}

\section{Introduction}

In \cite{alstrup} it is shown that there exists a graph $G$ with $O(n)$
vertices such that every forest on $n$ vertices is an induced subgraph of $G$.
In particular, then, every tree on at most $n$ vertices is an induced subgraph of $G$.  However, this $G$ is carefully constructed to have this property, and one might wonder whether there are graphs that naturally contain all large trees as induced subgraphs.  In this paper we show this for 
{\it prime sum graphs} \cite{chen} (where the ``large'' in ``all large trees" is somewhat smaller than for the \cite{alstrup} construction), which we will define later.  First, though, it's worth discussing how special the property of containing all large trees as induced subgraphs actually is by considering typical or random graphs:

In \cite{frieze} A. Frieze and B. Jackson proved that if
$p = c/n$, $c > 1$ a constant, then almost surely 
an Erd\H os-Re\'nyi random graph $G = G(n,p)$ contains an induced tree on at least 
$$
c^{-2} \min \{\sqrt{c-1}, 1/2\} n
$$
vertices.  

In light of the theorem on the so-called ``giant component" (see \cite{erdos}), if $c < 1$ then the largest connected component of $G$ has size $O(\log n)$, which immediately implies that the largest induced tree in $G$ has size at most $O(\log n)$.
\bigskip

We could ask a related question here:  suppose $T$ is a given tree on $m \leq n$ vertices.  Under what conditions on $m$ and $p$ does $G(n,p)$ contain an induced subgraph isomorphic to $T$ with probability $1-o(1)$?  

Consider the case where $T$ is the line graph on $m$ vertices (with edges $\{1,2\}, \{2,3\}, ..., \{m-1,m\}$ and no others).  Let us also assume $m > (2+\varepsilon)\log n$
(for reasons that will become clear later on).  For each choice of labels $x_1,...,x_m$ for an induced copy of $T$ in $G(n,p)$, the probability that
the edges $\{x_i,x_{i+1}\}$, $i=1,...,m-1$ are all in $G(n,p)$ is $p^{m-1}$; and the probability that none of the edges $\{x_i,x_j\}$, $|i-j| \geq 2$ are in $G(n,p)$ is $(1 - p)^{{m \choose 2} - (m-1)}\ =
(1-p)^{(m-1)(m - 2)/2}$. So, the expected number of 
induced copies of $T$ in $G(n,p)$ is at most (half, due to sequence reversals) the number of sequences $x_1,...,x_m$ times $p^{m-1} (1-p)^{(m-1)(m-2)/2}$, which is at most 
$$
n^m p^{m-1} (1-p)^{(m-1)(m-2)/2}\ = n^m p^{m-1} (1-p)^{m^2/2 - O(m)}.
$$
Now, if $p > (2+\varepsilon)(\log n)/m$, say, then this count will
be $o(1)$; and therefore a line graph on $m$ vertices is unlikely to be an induced subgraph of our $G(n,p)$.  Note that this is where we need the fact that $m > (2+\varepsilon) \log n$, which guarantees that $p$ doesn't exceed $1$.

So, it would be {\it highly improbable} that, say, $G(n,p)$ contains an induced subgraph isomorphic to a line graph even on just about $(2+\varepsilon)\log^2 n$ vertices when $p > (1+o(1))/\log n$.  Note that for $p$ this large, the expected degree of any vertex of $G(n,p)$ will by $(n-1)p > (1+o(1))n/\log n$.
\bigskip

Perhaps surprisingly, then, we establish that so-called {\it prime sum graphs} \cite{chen} on vertices $1,2,...,n$, which have average vertex degree $(1+o(1))n/\log n$, 
have the property that they contain an induced copy of each tree on $m$ vertices, for all $m \leq \exp(c (\log n)/\log\log n)$, for some $c > 0$.
\bigskip

\noindent {\bf Definition of the graph ${\cal P}_n$.}  Let $n \geq 1$, and consider the graph on vertex set $\{1,2,...,n\}$, where there is an edge connecting $i$ to $j$ if and only if $i+j$ is a prime number.  We will denote this graph by ${\cal P}_n$ and call it a ``prime graph on vertices $1,...,n$".
\bigskip

Before embarking on the discussion of the above result alluded to, it's worth mentioning that one shouldn't expect the average degree of a graph to tell you much about whether it contains induced copies of large trees, since one can imagine artificially-constructed graphs that are fairly dense, overall (the average vertex degree is near to $n$), that nonetheless contain a small, sparse subgraph where all the of the induced copies of those trees are found.  However, one {\it would} perhaps expect the average degree to tell you something in the case of random graphs (as we did above) or some type of quasi-random graph, as they have a more uniform structure.  Although prime graphs are not quasi-random, they do have a lot of uniform structure, after accounting for "local obstructions".  For example, from the Hardy-Littlewood Prime $k$-tuples Conjecture \cite{hardy}, given any integers $1 \leq i_1 < i_2 < \cdots < i_k$, with the property that they don't occupy all the congruence classes mod $p$ for any prime $p$, as $n$ tends to infinity there are $\sim \kappa(i_1,...,i_k) n (\log n)^{-k}$ integers $0 \leq i_0 < n - i_k$ so that $i_0i_1, ..., i_0i_k$ are all edges 
in the graph -- that is, $i_0 + i_1, ..., i_0 + i_k$ are all prime numbers.  Here $\kappa(i_1,...,i_k)$ is some constant depending only on $i_1,...,i_k$.  
\bigskip

It will turn out that our theorem regarding the structure of ${\cal P}_n$ does not require a lot of information about prime numbers, and that we can prove a similar theorem for much denser graphs, such as certain graphs ${\cal Q}_n(q)$ defined as follows:
\bigskip

\noindent {\bf Definition of the graph ${\cal Q}_n(q)$.}  Given an integer $q \geq 2$, we define ${\cal Q}_n(q)$ to be the graph on vertices $1,2,...,n$ where there is an edge connecting $i$ to $j$ if and only if $i+j$ is coprime to $q$.
\bigskip

Our main theorem of this paper alluded to above is as follows:
\bigskip

\begin{theorem}[Main Theorem]\label{maintheorem}  There exists an absolute constant $c > 0$ such that the following holds for all sufficiently large integers $n$: 
suppose $T$ is a tree on $m$ vertices where 
\begin{equation}\label{mbound}
m\ \leq\ \exp \left ( {c\log n \over \log\log n} \right ).
\end{equation}
Then, ${\cal P}_n$ and ${\cal Q}_n(q)$ both contain a copies of $T$ as induced subgraphs, where $q$ is the product of the primes in
$[5, (10 c + o(1)) \log n / \log(3/2)]$.  

An alternative way to say this is that: there exists an absolute constant $c' > 0$ so that for every $m \geq 1$, and any tree $T$ on $m$ vertices, we have that ${\cal P}_n$ and ${\cal Q}_n(q)$ contain copies of $T$ when
$$
n\ >\ \exp(c' m \log\log m).
$$
\end{theorem}

We do not believe that these bounds are anywhere near optimal, for the following reason:  the main information about the primes that the proof for ${\cal P}_n$ uses is the fact that a prime number $p > p_0$ is not divisible by any prime $q <10 c \log p /(\log 3/2)\log\log p$, and hence why the theorem gives the same results for
${\cal P}_n$ and ${\cal Q}_n$.  
\bigskip

We note that the prime number theorem implies that 
$$
q\ =\ \prod_{5 \leq p < 10c \log n / \log(3/2) \atop p\ {\rm prime}} p
\ =\ n^{10c/\log(3/2) + o(1)},
$$
and that the average vertex degree of ${\cal Q}_n(q)$ is
\begin{eqnarray*}
&& {1 \over n} \sum_{i=1}^n 
\#\{1 \leq j \leq n\ :\ {\rm gcd}(i + j,q)=1\}\nonumber \\
&&\ \ \ \ \gg\ {n \varphi(q)\over q}\ =\ 
n \prod_{5 \leq p < 10c \log n/\log(3/2) \atop p\ {\rm prime}}
\left ( 1 - {1 \over p} \right )\ \gg\ {n \over \log\log n}.
\end{eqnarray*}
(There is no point in being precise about the constants, since small modifications of the construction can yield slightly improved lower bounds on the average degree.)
\bigskip

We close the introduction with a two questions for further study:

\begin{enumerate}  

\item First, what is the true size of the largest $M = M(n)$ such that the prime graph ${\cal P}_n$ contains an induced copy of every tree on $m \leq M$ vertices?

\item Can one extend the main theorem to induced copies of graphs a little more complex than trees, where some cycles are added?  Obviously the subgraph can't contain an odd cycle, since ${\cal P}_n$ is bipartite (the even numbers up to $n$ form one part and the odd numbers up to $n$ 
form the other; if $a+b$ is an odd prime then one of $a$ or $b$ must be even, and the other must be odd).  We could, in fact, extend the definition of the graph to where $ab$ is an edge if $a+b$ is $p$ or $2p$ for some prime $p$, and then it is no longer bipartite, which would allow for more complex induced subgraphs.

\end{enumerate}

\section{Proof of the Main Theorem}

\subsection{Two Lemmas}

In order to prove the theorem we will need the following two lemmas.

\begin{lemma}\label{lemma1}
Suppose $T$ is a tree on $m \geq 3$ vertices, which we will also use to denote its vertex set (as we will all other trees mentioned).  Then, we can find two sub-trees $U_1$ and $U_2$ whose union of vertices is all those of $T$, such that $|U_1 \cap U_2| = 1$ (one vertex in common), and such that
$$
|U_1|,\ |U_2|\ \geq\ {m \over 3}.
$$
\end{lemma}

\noindent {\bf Remark:}  We note that this lemma is essentially best-possible, since in the case where $m=3k+1$ we can take $T$ to be the tree formed from three line or path graphs (graphs on vertices $v_1,...,v_k$ where $v_i$ connected to $v_j$ if and only if $|i-j|= 1$) of length $k$ and then connecting the ends of each to a final vertex $v$, which then has degree $3$.  This tree cannot be divided into two subtrees both of size larger than $m/3$ since:  two such disjoint trees would 
have at least $k+1$ vertices, which means both would have to include a vertex outside the paths of length $k$.  One tree would have to include the final vertex $v$ of degree $3$; but then once that vertex is used up, the other couldn't use it, and then you're stuck -- it couldn't have $k+1$ vertices if it doesn't contain $v$.

\begin{lemma}\label{lemma2}
Suppose $T$ is a tree on $m$ vertices.  It is possible to assign to every vertex $x$ in the tree $T$ a vector 
$$
v_x\ =\ (x_1,x_2,...,x_d)\ \in\ \{-1,1,2\}^d,
$$
where the following all hold:

\begin{enumerate}

\item $d\ \leq\ {10 \log m \over \log(3/2)}$

\item For each vertex $x$ in $T$ the number of times that
$1$ appears as a coordinate of $v_x$ is the same as $-1$ appears and that $2$ appears.  (In particular this means $3$ divides $d$.)

\item And, for every $y,z$ vertices of $T$,
$$
y\ {\rm adjancent\ to\ }z\ \Longleftrightarrow\ v_y+v_z \in 
\{-2,1,2,3,4\}^d.
$$
\end{enumerate}
\end{lemma}
\bigskip

\subsubsection{Proof of lemma \ref{lemma1}}

 We will form a sequence of vertices $v_1,v_2,v_3, ...$ and a sequence of trees $T_2,T_3, ...$ as follows:  we begin by letting $v_1$ be any leaf of $T$ and let $T_1 = T$.  Suppose we have constructed $v_1,...,v_k$.  We now show how to construct $v_{k+1}$ and the tree $T_{k+1}$.  Let $F$ denote the induced subgraph of $T$ gotten by deleting $v_k$.  We note that $F$ is a forest.  Let $T_{k+1}$ be the largest tree that is a subgraph of $F$, and let $v_{k+1}$ denote the unique vertex of $T_{k+1}$ that is connected to $v_k$ in $T$.  

Now, we can't have that the trees in the sequence $T_2, T_3, ...$ are decreasing in size indefinitely; and so, there exists $i$ such that 
\begin{equation}\label{Tii1}
|T_i| \leq |T_{i+1}|.
\end{equation}
Note that if this happens then $i \geq 2$, since:  $T_1 = T$, and $T_2$ is the tree gotten by removing $v_1$ from $T$ (recall $v_1$ is a leaf), together imply $|T_2| = |T_1|-1$.  So, $|T_i| > |T_{i+1}|$ for $i=1$.  

Another way to describe the trees $T_i$ and $T_{i+1}$ for this particular $i \geq 2$ where $|T_i| \leq |T_{i+1}|$ is as follows:  remove the edge connecting $v_{i-1}$ to $v_i$.  Then, the tree $T_i$ will be the part of that connected subgraph (tree) containing $v_i$, and $T_{i+1}$ will be the connected subgraph (tree) containing $v_{i-1}$.  So, 
\begin{equation}\label{Tisum}
|T_i| + |T_{i+1}| = m,
\end{equation}
and then together with (\ref{Tii1}) we deduce
$$
|T_i|\ \leq\ m/2.
$$

Now, if $m/3 \leq |T_i| \leq m/2$, then from (\ref{Tisum}) we deduce
$m/2 \leq |T_{i+1}| \leq 2m/3$; and so we can let 
$U_1 = T_{i+1}$ and then let $U_2$ be the tree gotten by adding to $T_i$ the vertex $v_{i-1}$ and edge from $v_{i-1}$ to $v_i$.

On the other hand, if $|T_i| < m/3$, then:  first, let $F'$ denote the forest
gotten by removing $v_{i-1}$ from $T$, and note that $T_i$ is the largest tree in $F'$.  Thus, all the trees in $F'$ have size $< m/3$ (not just $T_i$); and so, we can keep unioning one tree from $F'$ after another (unioning vertices and edges) until that union has size in the interval $[m/3, 2m/3)$; and then we can let $U_1$ be that collection of trees, together with $v_{i-1}$ and all its connections to them, producing a tree of size 
$|U_1| > m/3$.  We let $U_2$ be the union of the remainder of the trees -- where this union then has size at least $m/3$ also -- together with $v_{i-1}$ and its connections to those remaining trees.  Then, $|U_2| > m/3$, as well, and satisfies $|U_1 \cap U_2| = |\{v_{i-1}\}| = 1$; and every vertex of $T$ is either in $U_1$ or $U_2$.  This completes the proof.  
\hfill $\blacksquare$

\subsubsection{Proof of lemma \ref{lemma2}} \label{lemma2section}

In the case where the tree $T$ has $3$ vertices, we choose $d=3$ and take 
$v_x,v_y, v_z$ to be the following
$$
v_x\ =\ (-1,1,2),\ v_y\ =\ (2,1,-1),\ v_z\ =\ (1,2,-1).
$$
Note that $v_x + v_y$ and $v_y + v_z$ are not $0$ in any coordinate, yet
$v_x + v_z$ is $0$ in the first coordinate.  We also solve the case when $T$
has $2$ vertices, $x$ and $y$, using $d=3$ and keeping the same $v_x, v_y$.

Assume, for proof by induction, we've proved the {\bf Claim} above for all trees $T$ on $2 \leq m \leq k$ vertices, $k \geq 3$.  Now we prove it for a tree on $m=k+1 \geq 4$ vertices.  

We begin by applying Lemma \ref{lemma1} to our tree $T$.  We therefore have that $T$ can be expressed as the union of two sub-trees $U_1$, $U_2$, with $|U_1 \cap U_2| = 1$, where for $i=1,2$,
$$
{m \over 3}\ \leq\ |U_i|\ \leq\ {2m\over 3}+1.
$$
If $m=k+1=4$ we can do better:  the lower bound on $|U_i|$ from
Lemma \ref{lemma1} implies $|U_1|, |U_2| \geq 2$; and therefore since $|U_1| + |U_2| = m+1 = 5$, we also would have $|U_1|, |U_2| \leq 3 \leq k$.  So, we could apply the induction hypothesis to $U_1$ and $U_2$.

Also, if $m \geq 5$ we have that $\lfloor 1 + 2m/3\rfloor \leq m-1$, and so
we may apply the induction hypothesis in that case as well.  

So, let us assume we have an assignment of vectors $v_x$ of dimension $d_1$ to all vertices $x \in U_1$, and an assignment of vectors $w_y$ of dimension $d_2$ to all the vertices $y \in U_2$.  Now, we may have that $d_1 \neq d_2$.  Without loss of generality, let's assume $d_1 \leq d_2$.  

If $d_1 < d_2$, then we can pad the vectors $v_x$ with the pattern of coordinates $-1,1,2$ again and again as needed (recall that $3$ divides $d_1$ and $d_2$), until we produce vectors of length $d_2$.  For example, suppose
$$
v_x = (x_1,x_2,x_3,...,x_{d_1}),
$$
and that $d_2 = d_1+6$.  Then, after padding, we get the new vector
$$
(x_1,x_2,...,x_{d_1},-1,1,2,-1,1,2).
$$
And note that after this padding is added, we still get condition 3 of Lemma \ref{lemma2} holding for the tree $U_1$.

Let $u$ denote the vertex in common between $U_1$ and $U_2$.
We note that $v_u$ and $w_u$ are, in general, not equal.  However, they both have exactly $d_2/3$ coordinates equal to $-1$, $d_2/3$ coordinates equal to $1$, and $d_2/3$ coordinates equal to $2$.  Thus, there is a permutation $\sigma$ of the coordinates of $v_u$ so that $\sigma(v_u) = w_u$.  

We note that if we define the new vector mapping 
$$
v'_x\ =\ \sigma(v_x)
$$
gotten by permuting the coordinates of $v_x$, then $v'_x$ satisfies all the properties of Lemma \ref{lemma2} for the tree $U_1$ when used in place of $v_x$.

Now we define a new set of vectors $\rho_x$ for $x$ a vertex of $T$ as follows:  if $x \in U_1$ and $x \neq u$, and if write
$$
v'_x\ =\ (x_1, ..., x_{d_2}),
$$
then we let
$$
\rho_x\ :=\ (x_1,...,x_{d_2},-1,1,2).
$$
Next, if we write
$$
v'_u\ =\ w_u\ =\ (u_1,...,u_{d_2}),
$$
then we let
$$
\rho_u\ :=\ (u_1,...,u_{d_2},2,1,-1).
$$
Finally, if $y\in U_2$, $y \neq u$, and 
$$
w_y\ =\ (y_1,...,y_{d_2}),
$$
then we let 
$$
\rho_y\ :=\ (y_1,...,y_{d_2},1,2,-1).
$$

We claim that these definitions will ensure that $\rho$ gives an assignment of vectors to the vertices of $T$ that satisfy the conclusion of Lemma \ref{lemma2}.  We first note that if $x\in U_1$, $y\in U_2$, with $x,y \neq u$, then we don't have $x$ adjacent to $y$; and, indeed, when $\rho_x + \rho_y$ is $0$ in the third-to-last coordinate.  Also, if $x\in U_1$ is adjacent to $u$, then none of the digits of $v_x + v_u$ are $0$, by the induction hypothesis; and the same will be true of $v'_x + v'_u$, and therefore also $\rho_x + \rho_u$.  The same holds for $\rho_y + \rho_u$ if $y \in U_2$.  Thus, the third condition of Lemma \ref{lemma2} is satisfied.  

It remains to bound the dimension or number of coordinates of $\rho_x$ for $x$ a vertex of $T$:  we clearly have that the dimension $d$ of the vectors $\rho_x$ satisfies
$$
d\ =\ \max \{d_1,d_2\} + 3.
$$
The $3$ comes from the three extra digits of padding in transforming $v'_x$ into $\rho_x$ or $w_y$ into $\rho_y$.  By the induction hypothesis,
\begin{eqnarray*}
d\ \leq\ \max \left \{ {10\log |U_1| \over \log(3/2)},\ {10\log |U_2| \over \log(3/2)}\right \} + 3
\ &\leq&\ {10\log(2m/3+1) \over \log(3/2)} + 3\\
&=&\ {10\log (m+ 3/2) \over \log(3/2)}-7.
\end{eqnarray*}
We wish to see for which $m$ this is at most $10\log(m)/\log(3/2)$; that is, we seek $m$ where
$$
10\log(m+3/2) - 10 \log m\ <\ 7 \log(3/2).
$$
That is, we seek $n$ where
$$
\log(1 + 3/2m)\ <\ {7 \log(3/2) \over 10}\ =\ 0.12326...
$$
When $m=5$ we already have $\log(1+3/2m) = 0.1139...$, and larger values
of $m$ will only make the left-hand-side of the inequality smaller, so
it holds for all $m \geq 5$.  When $m=4$ we saw earlier
in the proof of this lemma that $|U_1|, |U_2| \leq 3$; and so in that case we have
$$
d\ \leq\ {10 \log 3 \over \log(3/2)} + 3\ <\ 30.095...\ <\ 
34.190...\ =\ {10 \log 4 \over \log(3/2)}.
$$
The bound $d < 10 \log(m)/\log(3/2)$ also holds for $m=2$ (in fact, we can use
$d=1$ in that case).  So, the induction step is proved.

\hfill $\blacksquare$

\subsection{Conclusion of the proof of the Main Theorem}

We will need the following estimate due to Fogels \cite[p. 83]{fogels}, which gives interval lower-bound estimates for primes in arithmetic progressions, which can be thought of as an extension of Linnik's Theorem:

\begin{theorem}[Fogels's bound]\label{fogelstheorem} There exist constants 
$0 < c_2 < c_1 < 1 < c_3$ such that the following holds for all integers 
$q \geq 2$ and all $1 \leq a \leq q-1$ with gcd$(a,q)=1$:  if $x > q^{c_3}$, then
$$
\pi([x, x+x^{c_1}]; q,a)\ >\ x^{c_2},
$$
where $\pi(I;q,a)$ denotes the number of primes in an interval $I$ that are $a \pmod{q}$.
\end{theorem}

\noindent {\bf Remark:}  We really only use the fact that $\pi([x,x+x^{c_1}]; q,a) > 0$, which seems like a traditional Linnik's Theorem result, except that we have restricted to an interval that can be a small power of $q$ away from $a$.
\bigskip

Next, we apply Lemma \ref{lemma2} to the tree $T$.  Since
the number of vertices in $T$ is 
$$
m\ \leq\ \exp \left ( {c\log n \over \log\log n} \right ),
$$
where $c > 0$ will be determined later and will depend
on $c_1,c_2,c_3$ in Theorem \ref{fogelstheorem}.  

We assume that the vertices of $T$ are labeled $1,2,3,...,m$, and we order them so that the vertex labeled $j \geq 2$ is connected to some vertex labeled $i$ where 
$1 \leq i < j$.  For a vertex labeled $h$ in $T$, we write the vector encoding 
$\rho_h$ given by the lemma:
\begin{equation}\label{rhoh}
\rho_h\ =\ (h_1,...,h_d).
\end{equation}
Note that
$$
d\ \leq\ {10 \log m \over \log(3/2)}\ <\ 
{10c \log n \over \log(3/2) \log\log n}.
$$

Now, let $q$ be the product of the $d$ consecutive prime numbers starting at $5$.  So, for example, if $d=3$, then $q = 5 \cdot 7 \cdot 11$.  As a consequence of the prime number theorem, we have that the product of the first $k$ primes
has size $\exp((1+o(1))k \log k)$, and therefore
\begin{equation}\label{qnbound}
q\ <\ \exp((1+o(1))d \log d)\ <\ 
n^{(1+o(1)) 10c/\log(3/2)}.
\end{equation}

Let $q_1 < q_2 < \cdots < q_d$ denote those primes making
up $q$.

Using the Chinese Remainder Theorem we let $a_h$ denote
the unique integer satisfying
$$
1 \leq a_h \leq q-1,\ {\rm and\ for\ all\ }i=1,...,d,\ 
a_h\ \equiv\ h_i \pmod{q_i},
$$
where $h_i$ is as in (\ref{rhoh}).
We note since the $h_i = -1,1$ or $2$, and $q_i \geq 5$, that ${\rm gcd}(a_h, q) = 1$.  Also note that, as a consequence of Lemma \ref{lemma2} part 3, 
if $x$ and $y$ are labels of vertices in $T$ that are connected, then 
$a_x + a_y$ is coprime to $q$, since for each $i=1,2,...,d$ 
this sum will be congruent to 
$-2, 1, 2, 3,$ or $4$ mod $q_i$, none of which is $0$; and if $x$ and $y$ are {\it not}
connected in $T$, then $a_x + a_y$ has a non-trivial common divisor with $q$ -- that is, for some $i$ we will have $a_x + a_y \equiv 0 \pmod{q_i}$.

We now choose distinct integers $j_1,...,j_m$ (that we think of as labels of vertices in the prime graph), one for each vertex of $T$.  First, choose $j_1 = a_1 + q\lceil q^{c_3-1}\rceil$ (this term $q\lceil q^{c_3-1}\rceil$ guarantees that $j_1 > q^{c_3}$ while being congruent to $a_1 \pmod{q}$).  Then, suppose we have chosen $j_1, j_2, ..., j_k$, $k \leq m-1$ where
\bigskip

\noindent {\bf Induction Hypothesis:}
\begin{enumerate}

\item Each $j_i \equiv a_i \pmod{q}$.

\item Each $j_i$ satisfies $q^{c_3} \leq j_i < 2 q^{c_3}$.

\item All $j_1,j_2, ..., j_k$ are distinct.

\item And finally, if $i, i' \leq m-1$ are adjacent in 
$T$, then $j_i + j_{i'}$ is a prime number, meaning that $j_i$ and $j_{i'}$ are
adjacent in prime graph.  Note also that if $i$ and $i'$ are {\it not} adjacent in $T$, then $a_i + a_{i'}$ has a non-trivial common divisor with $q$, making $j_i + j_{i'}$ (which is congruent to $a_i + a_{i'} \pmod{q}$) not a prime number, so that there is no edge connecting $j_i$ to $j_{i'}$ in the prime graph.  

\end{enumerate}
\bigskip

Now we show how to choose $j_{k+1}$:  by our labeling scheme for vertices of $T$, the vertex $k+1$ in $T$ is connected to some vertex $i \leq k$.  It follows, then, that $a_{k+1} + a_i$ is coprime to $q$.  Thus, the arithmetic progression $a_{k+1} + a_i + t q$, $t=0,1,2,...$ contains infinitely many prime numbers.  We now
apply Theorem \ref{fogelstheorem} using $x=j_i + q^{c_3}$ and
$a\equiv a_{k+1} + a_i\pmod{q}$, along with our induction hypothesis 
(in particular that $j_i > q^{c_3}$) and deduce the existence of a prime
number  
\begin{equation}\label{Pinterval}
P\ \in\ [j_i + q^{c_3},\ j_i + q^{c_3} + (j_i + q^{c_3})^{c_1}]
\end{equation}
such that $P \equiv a \equiv a_{k+1} + a_i \pmod{q}$.
Note that $P$ has the form $a_{k+1} + a_i + t q$, $t \in {\mathbb Z}$.  
We then let 
$$
j_{k+1} = P -  j_i\ \neq\ j_i.
$$
($P - j_i \neq j_i$ because otherwise $P = 2j_i$, which is not prime.)
Note that this implies
$$
j_{k+1}\ \equiv\ a - j_i\ \equiv\ a - a_i\ \equiv\ a_{k+1} \pmod{q},
$$
which is the first condition we need to show for the induction step of our induction
proof.  

Next, we also note that since $P \geq j_i + q^{c_3}$, it follows that
$j_{k+1}\ \geq\ q^{c_3}$, which is part of the second part of the induction step
we need to prove (see part 2 of the Induction Hypothesis above).  Also,
from (\ref{Pinterval}) and part 2 of the induction hypothesis,
$$
j_{k+1}\ =\ P- j_i\ \leq\ q^{c_3} + (j_i + q^{c_3})^{c_1}\ \leq\ 
q^{c_3} + (3q^{c_3})^{c_1}\ <\ 2 q^{c_3}
$$
for $q$ sufficiently large.

We have that $j_{k+1} + j_i = P$, so that $j_{k+1}$, $j_i$ are connected in the prime graph, and the same is true of sums
$j_i + j_{i'}$ where $i$ and $i'$ are connected in $T$, by the 
induction hypothesis.  Thus, the fourth property of the induction step is proved.  

Finally, we establish the third property:  it suffices to show
that $j_{k+1}$ is distinct from $j_1,...,j_k$, since by the induction hypothesis $j_1,...,j_k$ are all distinct already.  
Note once again that 
$j_{k+1} \neq j_i$ (recall we are assuming $i$ and $k+1$ are connected in $T$) since otherwise $2 | P$.  
If $i' \leq k$, $i' \neq i$, then $k+1$ and $i'$ are not connected in $T$, so $a_{k+1} + a_{i'}$ has a non-trivial common divisor with $q$.  We have then that $j_{k+1} \neq j_{i'}$ since otherwise we would have 
$$
2 j_{k+1}\ \equiv\ j_{k+1} + j_{i'}\ \equiv a_{k+1} + a_{i'} \pmod{q}
$$
has a common factor with $q$, yet $2j_{k+1} \equiv 
2\ \cdot \{-1,1,2\} \pmod{q'}$ for
any $q' | q$, so can't be $0$ mod $q'$ for any such $q'$.  
Thus, we can conclude that all the $j_i$ are distinct, which is the third condition we needed to show the induction step. 
\bigskip

The only thing that remains to be shown to prove the theorem is that all the 
$j_i \leq 2 q^{c_3} \leq n$, so that they are vertices of ${\cal P}_n$ and ${\cal Q}_n(q)$.  From (\ref{qnbound}) we see this will hold for any 
$0 < c < \log(3/2)/10 c_3$.  
$\hfill$ $\blacksquare$

\end{document}